\documentclass[11pt,leqno,fleqn,epsfig]{article}
\usepackage{amssymb,epsfig, amsmath}
\usepackage{mathrsfs}  

\newfont{\graf}{eufm10}

\newcommand{\R}{\mathbb{R}}
\newcommand{\N}{\mathbb{N}}

\newcommand{\ba}{\boldsymbol a}
\newcommand{\bb}{\boldsymbol b}
\newcommand{\bc}{\boldsymbol c}
\newcommand{\bd}{\boldsymbol d}
\newcommand{\bfe}{\boldsymbol e}
\newcommand{\bbf}{\boldsymbol f}
\newcommand{\bg}{\boldsymbol g}
\newcommand{\bh}{\boldsymbol h}
\newcommand{\bi}{\boldsymbol i}
\newcommand{\bj}{\boldsymbol j}
\newcommand{\bk}{\boldsymbol k}
\newcommand{\bl}{\boldsymbol l}
\newcommand{\bm}{\boldsymbol m}
\newcommand{\bn}{\boldsymbol n}
\newcommand{\bo}{\boldsymbol o}
\newcommand{\bp}{\boldsymbol p}
\newcommand{\bq}{\boldsymbol q}
\newcommand{\br}{\boldsymbol r}
\newcommand{\bs}{\boldsymbol s}
\newcommand{\bt}{\boldsymbol t}
\newcommand{\bu}{\boldsymbol u}
\newcommand{\bv}{\boldsymbol v}
\newcommand{\bw}{\boldsymbol w} 
\newcommand{\bx}{\boldsymbol x}
\newcommand{\by}{\boldsymbol y}
\newcommand{\bz}{\boldsymbol z}
\newcommand{\buu}{  \underline{\boldsymbol u}}

\newcommand{\bA}{\boldsymbol A}
\newcommand{\bB}{\boldsymbol B}
\newcommand{\bC}{\boldsymbol C}
\newcommand{\bD}{\boldsymbol D}
\newcommand{\bE}{\boldsymbol E}
\newcommand{\bF}{\boldsymbol F}
\newcommand{\bG}{\boldsymbol G}
\newcommand{\bH}{\boldsymbol H}
\newcommand{\bI}{\boldsymbol I}
\newcommand{\bJ}{\boldsymbol J}
\newcommand{\bK}{\boldsymbol K}
\newcommand{\bL}{\boldsymbol L}
\newcommand{\bM}{\boldsymbol M}
\newcommand{\bN}{\boldsymbol N}
\newcommand{\bO}{\boldsymbol O}
\newcommand{\bP}{\boldsymbol P}
\newcommand{\bQ}{\boldsymbol Q}
\newcommand{\bR}{\boldsymbol R}
\newcommand{\bS}{\boldsymbol S}
\newcommand{\bT}{\boldsymbol T}
\newcommand{\bU}{\boldsymbol U}
\newcommand{\bV}{\boldsymbol V}
\newcommand{\bW}{\boldsymbol W}
\newcommand{\bX}{\boldsymbol X}
\newcommand{\bY}{\boldsymbol Y}
\newcommand{\bZ}{\boldsymbol Z}

\newcommand{\bPhi}{\boldsymbol \Phi}

\newcommand{\mes}{\operatorname{\rm mes}}
\newcommand{\esssup}{\operatorname*{ess\,sup}}

\newcommand{\diver}{\operatorname*{div}}
\newcommand{\supp}{\operatorname*{supp}}

\newcommand{\const}{\operatorname*{const}}

\newcommand{\Beweisende}{\rule{0.2cm}{0.2cm}}

\newcommand{\be}{\begin{equation}}
\newcommand{\ee}{\end{equation}}
\newcommand{\bea}{\begin{eqnarray}}
\newcommand{\eea}{\end{eqnarray}}
\newcommand{\bean}{\begin{eqnarray*}}
\newcommand{\eean}{\end{eqnarray*}}

\newcommand{\intl}{\int\limits}

\newcounter{secnum}

\begin{document}

\title{A note on the pressure of strong solutions to the Stokes system  
\\
in bounded and exterior domains}
\date{}
\author{J\"org Wolf\\
Department of Mathematics \\
Humboldt University Berlin\\
joerg@math.hu-berlin.de}

\maketitle

\begin{abstract}
We consider the Stokes problem in an exterior domain $\Omega \subset \R^n$ 
with an external force $\bbf \in L^s(0,T; \bW^{k,\, r}(\Omega ))\, 
(k\in \N, 1<r<\infty)$. 
 In the present paper we show that in contrast to $\bu$ the boundary regularity of the pressure can be improved according to 
 the differentiability  of  $\bbf$ up to order $ k$.  In particular, this implies that the pressure is smooth with respect to $x\in \Omega$  
if $\bbf$ is smooth with respect to $x\in \Omega $.
\end{abstract}

\vspace{0.5cm}

\vspace{0.5cm}

\vspace{0.5cm}
\section{Introduction}
\label{sec:1}
\setcounter{secnum}{\value{section} \setcounter{equation}{0}
\renewcommand{\theequation}{\mbox{\arabic{secnum}.\arabic{equation}}}}

Let $\Omega \subset  \R^n\, (n\in \N, n\ge 2)$  be an exterior domain, i.\,e.  
$ \R^n\setminus \overline{\Omega }$ is a bounded domain in $\R^n$.  
Let $0<T< +\infty$. Set 
$Q= \Omega \times (0,T)$.  In the present paper we consider the Stokes problem 
\begin{align}
  \diver \bu &=0  \quad  \mbox{in}\quad Q
\label{1.1}
\\
\partial_t \bu - \Delta \bu  &=  - \nabla p  + \bbf\quad \mbox{in}\quad Q,
\label{1.2}
\\
\bu  &=  0 \quad \mbox{on}\quad \partial \Omega \times (0,T),\quad 
\lim_{|x|\to \infty} \bu (x, \cdot)  =  0,
\label{1.3}
\\
\bu(0)  &=  0\quad\mbox{in}\quad \Omega,
\label{1.4}
\end{align}
where $\bu=(u^1,\ldots , u^n)$ denotes the unknown velocity of the fluid, $p$ the unknown pressure and $\bbf$ the  given external force.  The Stokes problem 
has been extensively studied in the past.  In particular, for  the case $\Omega$ 
is the half space or an $C^2$  domain with compact boundary the $L^p$-theory 
is well-known. Based on potential theory in \cite{SOL}  Solonnikov proved that for every $\bbf\in \bL^q(Q)$ there exists a unique solution 
$(\bu, p) $  to \eqref{1.1}--\eqref{1.4} such that  
$\partial_t \bu, \nabla^2 \bu \in \bL^q(Q)$,  and  $\nabla p \in \bL^p(Q)$.  
By using the semi group approach, similar results have been obtained in 
\cite{GI}, \cite{GI2},  \cite{BS}.  For the corresponding estimates on the pressure we 
refer to \cite{SOVW}.   An optimal result for the anisotropic case when 
$\bbf$ belongs to $ L^s(0,T; \bL^q(\Omega ))$  has been proved in \cite{GS} for 
the cases $\Omega =\R^n$, $\Omega =\R^n_+$,  and   a $ C^2$ domain  $\Omega $ with  compact boundary.  

\hspace*{0.5cm}
By standard arguments from the regularity theory of parabolic equations one gets the 
regularity $\bu$ and $p$  in dependence of the  regularity of the right-hand side $\bbf$ in time and space.  
However, if $\bbf$ is only smooth in $x\in \Omega $ it is not clear whether $\bu$ is smooth in $x$ up to the boundary.  In the present paper we will see that such a property at least holds for the 
pressure $p$, which is due to the fact that $\Delta p=0$ if $\diver \bbf =0$. More precisely, the condition$\bbf \in L^s(0,T; \bW^{k,\, q}(\Omega ))$\, 
$(1<s,q, <+\infty; k\in \N)$  implies $\nabla p \in L^s(0,T; \bW^{k,\, q}(\Omega ))$.    
Note that our result relies essentially on the fact that the initial data is zero.  
In general our result may not be true  as there is a counter-example 
obtained in \cite{HW}. More precisely,  there exists an initial data,  and 
a solution $\bu, p$ to the Stokes system 
such that  $\|\nabla \bu(t)\|_{\bL^2} $  is continuous as  $t \rightarrow 0^+$, while 
the corresponding estimate on the pressure $ \|p(t)\|_{L^2} $  may blow up as $t \rightarrow 0^+$.

\vspace{0.2cm}
\hspace*{0.5cm}
First we shall introduce the basic notations regarding the  function spaces  used throughout   the paper. By $W^{k,\, q}(\Omega ), W^{k,\, q}_0(\Omega )$ we denote 
the usual Sobolev spaces. Vector functions and spaces of vector valued functions will be denoted by bold face letters, i.\,e. we write 
$\bL^q(\Omega ), \bW^{k,\, q}(\Omega ), $ etc.   instead of $L^q(\Omega ; \R^n), 
W^{k,\, q}(\Omega ; \R^n)$, etc.   In addition, we use the following spaces 
of solenoidal functions 
\begin{align*}
\bL^q_\sigma (\Omega )  &= \,\, 
\mbox{closure of $\mathscr{C}^\infty_{0,\sigma } (\Omega )$ w.r.t. the norm $\|\cdot\|_{L^q} $} 
\\
\bW^{k,\, q}_{0,\sigma}  (\Omega ) & = \,\, \mbox{closure of $\mathscr{C}^\infty_{0,\sigma } (\Omega )$ w.r.t. the norm $\|\cdot\|_{W^{k, q} } $},
\end{align*} 
where $\mathscr{C}^\infty_{0,\sigma } (\Omega )$ stands for the space of all smooth 
solenoidal vector fields with compact support in $\Omega $.  Given a Banach space 
$X$ by $L^q(0,T; X)$ we denote the space of Bochner measurable functions 
$f: (0,T) \rightarrow X $ such that  
\begin{align*}
\|f\|_{L^q(0,T;X)}  ^q &= \intl_{0}^{T} \|f(t)\|_X^q d t < +\infty \quad  \mbox{if}\quad 1\le q <+\infty,
\\
\|f\|_{L^\infty(0,T;X)}  &= \esssup_{t\in (0,T)} \|f(t)\|_X < +\infty \quad  \mbox{if}\quad q =+\infty.    
\end{align*}

Now, let us introduce the notion of a strong solution to \eqref{1.1}--\eqref{1.4}.

\vspace{0.5cm}
{\bf Definition\,1.1}  Let $\bbf\in L^s(0,T; \bL^q(\Omega ))\, (1< s,q<+\infty)$.  
A pair $(\bu, p)$ is called a {\it strong solution} to \eqref{1.1}--\eqref{1.4} 
if $\bu\in L^s(0,T; \bW^{1,\, q}_{0, \sigma } (\Omega )), 
p \in L^s(0,T; L^1_{\rm loc} (\overline{\Omega}))$ and 
\[
\partial_i\partial_j\bu, \partial _t\bu,  \nabla p \in 
L^s(0,T; \bL^{ q}(\Omega )),\quad  i,j=1,2,3,
\]
such that \eqref{1.1}, \eqref{1.2} holds a.\,e. in $Q$, while 
\eqref{1.4} is fulfilled  such that  $\bu=0$ a.\,e. in $\Omega \times \{0\}$. 

\vspace{0.5cm}
For the existence of a strong solution to \eqref{1.1}--\eqref{1.4} cf. in \cite{GS}.

\vspace{0.5cm}
\hspace*{0.5cm}
Our main result is the following 

\vspace{0.3cm}
{\bf Theorem\,1}\, {\it Let $\Omega \subset \R^3$ be an exterior domain or a bounded domain with 
$\partial \Omega \in C^{2+k} \, (k\in \N)$.  For  $\bbf \in 
L^s(0,T ; \bW^{k,\, q}(\Omega ))$ $(1< s,q< +\infty; k\in \N)$, 
let $(\bu, p)$ be the strong solution to \eqref{1.1}--\eqref{1.4}. 
Then, 
\[
\nabla p \in L^s(0,T ; \bW^{k,\, q}(\Omega )).
\]
In addition, there holds
\be\label{1.5}
\|\nabla^{k+1}  p\|_{L^s(0,T; \bL^q(\Omega ))}   \leq  c
\|\bbf\|_{L^s(0,T;  \bW^{k,\, q}(\Omega ))},
\ee
where $c=\const>0$ depending only on $s,q, k$ and the geometric properties  of 
$\partial \Omega $.}

\section{Remarks on the equation $\diver \bv = f$}
\label{sec:2}
\setcounter{secnum}{\value{section} \setcounter{equation}{0}
\renewcommand{\theequation}{\mbox{\arabic{secnum}.\arabic{equation}}}}

Let $G\subset \R^n$  be a bounded domain, star-shaped with respect to a ball $B_R$.  
It is well known that for all $f\in L^q(G)$ with $(f)_G=0$ \,\footnotemark\, 
the equation $\diver \bv =f$  has a solution $\bv\in \bW^{ 1,\, q}_0(G)$ such that 
\footnotetext{\,Let $A\subset \R^n$ be a measurable set with $\mes(A)$. 
Given $v\in L^1(A)$ by $(v)_A$ we denote the mean value $\frac {1} {\mes(A)}
\intl_{A} v(x)dx.$}
\[
\|\nabla \bv\|_{L^q(G)} \le c \|f\|_{L^q(G)}
\]
with $c=\const>0$, depending on $n, q$ and $G$ (cf. \cite{BO}, \cite{GA}).   
In fact, the constant $c$ depends on the geometric property of $G$, namely the ratio of $G$ which is defined by 
\[
{\rm ratio}(G) := \frac {R_a(G)} {R_i(G)}, 
\]
 where
\begin{align*}
R_a(G) &=\inf \{R>0 \,|\, \exists B_R(x_0): G \subset B_R(x_0)\},  
\\
R_i(G) &=\sup \{r>0 \,|\, \exists B_r(x_0): G\,\, \mbox{is star-shaped w.r.t $B_r(x_0)$}\}.
\end{align*}
For instance ${\rm ratio}(G)=1$ if $G$ is a ball,  and ${\rm ratio}(G)=\,\sqrt[]{n}$ if 
$G$ is a cube.  Moreover, the ratio is invariant under translation and scaling, i.\,e.
\[
{\rm ratio}(\lambda G) ={\rm ratio}(G) \quad \forall\,  \lambda >0.
\]

\hspace*{0.5cm} 
Now, let $G$  such that $2< R_i(G)< 3$. In particular, $G$ is star shaped with respect to a ball 
$B_2=B_2(x_0)$.  Without loss of generality we may assume that $x_0=0$.  
Let $\phi \in C^\infty_0(B_2)$. We define 
\[
\mathscr{B}_{\phi }  f (x) =\intl_{\R^n} f(x-y) \bK_\phi (x,y) dy,\quad 
x\in \R^n,\quad f\in C^\infty_0(G),  
\]
where 
\[
\bK_\phi (x,y) = \frac {y} {|y|^{n} } 
 \intl_{0}^{\infty} \phi \Big(x+ r\frac {y} {|y|^n}\Big) (|y|+r)^{n-1}   d r,\quad 
(x,y)\in \R^n\times (\R^n\setminus \{0\}).
\]
As in \cite{BO}, \cite{GA}  it has been proves that $\mathscr{B}_{\phi }  f \in \bC^\infty_0(G)$  for all $f\in C^\infty_0(G)$. In addition, there holds
\be\label{2.2}
\|\nabla^k \mathscr{B}_{\phi }  f \|_{\bL^q(G)} \le c \|\nabla^{k-1} f\|_{L^q(G)} \quad \forall\, f\in C^\infty_0(G)
\ee
with a constant depending on $n, k, q, \phi $ and ${\rm ratio}(G)$ only.  Furthermore, 
there holds 
\be\label{2.3}
\diver \mathscr{B}_{\phi }  f = f \intl_{B_1} \phi(y)  dy  - \phi 
\intl_{G} f(y)  dy \quad \mbox{in} \quad G.
\ee
In particular, if $\intl_{B_1} \phi(y)  dy =1 $ and $\intl_{G} f(y)  dy =0$ then 
$\bv=\mathscr{B}_{\phi }  f $  solves the equation $\diver \bv=f$.   
Finally, by \eqref{2.2}  we may extend $\mathscr{B}_{\phi }$ to an operator 
$\mathscr{L}(W^{k-1, q}(G), \bW^{k,\, q}(G))$ denoted again by $\mathscr{B}_{\phi }$.  

\hspace*{0.5cm}
Let $i,j \in \{1, \ldots, n\}$.  Observing, that 
\begin{align*}
\partial_j\mathscr{B}_{\phi }  (\partial_i f)  &=
\partial_i\partial_j  \mathscr{B}_{\phi }  (f) - \partial_j\mathscr{B}_{\partial_i\phi }  (f) \quad \mbox{in} \quad G,\
\\
\partial_i\partial_j  \mathscr{B}_{\phi }  (f)  &=\partial_i\mathscr{B}_{\phi }  (\partial_j f)  + \partial_i\mathscr{B}_{\partial_j\phi }  (f) \quad \mbox{in} 
\quad G,
\end{align*}
we see that 
\[
\partial_j\mathscr{B}_{\phi }  (\partial_i f)  =
\partial_i\mathscr{B}_{\phi }  (\partial_j f)  + \partial_i\mathscr{B}_{\partial_j\phi }  (f) - \partial_j\mathscr{B}_{\partial_i\phi }  (f) \quad \mbox{in} \quad G.
\]
By the aid of  \eqref{2.2}, and Poincar\'{e}'s inequality, using the above identity, we get 
\begin{align}
&\| \partial_j\mathscr{B}_{\phi }  (\partial_i f)\|_{\bL^q(G)} \le 
c (\| \partial_j f\|_{L^q(G)} +\|f\|_{L^q(G)})\quad \forall\,  
f\in W^{ 1,\, q}_0(G),  
\label{2.4}
\\[0.3cm]
& \hspace*{-0.3cm} \begin{cases}
\|\nabla^2  \partial_j\mathscr{B}_{\phi}  (\partial_i f)\|_{\bL^q(G)} 
 \leq  
c (\| \partial_i \nabla_\ast \nabla    f\|_{\bL^q(G)}+ \|\partial_j  \partial_n 
\partial_n f\|_{L^q(G)}   +  
\|\nabla^2 f\|_{\bL^q(G)}) \,\footnotemark\,
\\[0.3cm]
\forall\,  f\in W^{3,\, q}_0(G),
\end{cases}
\label{2.4-1}
\end{align}
\footnotetext{\,Here $\nabla_\ast$ denotes the reduced gradient  
$(\partial_1, \ldots, \partial_{n-1} )$.}
where $c=\const>0$,  depending on $n, q$  and ${\rm ratio}(G)$.

\hspace*{0.5cm}
Now, let $G$ be a bounded domain, star-shaped with respect to a ball $B$. 
Let $R:= \frac {1} {2} R_i(G)$.  Thus, there exist $B_R(x_0)$ such that $G$ is star shaped to the ball $B_R(x_0)$. Without loss of generality we may assume that $x_0 =0$.  
Let $\phi \in C^\infty_0(B_1)$  with $\intl_{B_1} \phi (y) dy =1 $.  
We define 
$\mathscr{B} : W^{ k-1,\, q}_0(G) \rightarrow \bW^{ k,\, q}_0(G)$  by setting 
\[
\mathscr{B} (f)(x) =R \mathscr{B}_\phi  (\widetilde{f } ) \Big(\frac {x} {R}\Big),\quad x\in G,\quad  f \in W^{ k-1,\, q}_0(G),
\]
where $\widetilde{f }(y) = f(Ry)\, (y\in R^{-1} G)$. Using the transformation formula of the Lebesgue integral, in view of \eqref{2.2},  we see that 
\begin{align}
\|\nabla^k\mathscr{B} (f)\|_{\bL^q(G)} &= 
R^{n/q-k+1} 
\|\nabla^k\mathscr{B}_\phi  (\widetilde{f })\|_{\bL^q(R^{-1} G)}
 \leq  c R^{n/q-k+1} \|\nabla^{k-1} \widetilde{f }\|_{\bL^q(R^{-1} G)}
\cr   
 &=  c\|\nabla^{k-1} f \|_{\bL^q(G)},
\label{2.5}
\end{align}
where $c=\const>0$  depends on $n,q $ and ${\rm ratio}(R^{-1}G ) = {\rm ratio}(G)$.  
In addition,  from \eqref{2.4}, and \eqref{2.4-1} we deduce 
\begin{align}
&\| \partial_j\mathscr{B}  (\partial_i f)\|_{\bL^q(G)} \le  
c (\| \partial_j f\|_{L^q(G)} + R^{-1} \|f\|_{L^q(G)})\quad 
\forall\,  f\in W^{ 1,\, q}_0(G),  
\label{2.6}
\\[0.3cm]
& \hspace*{-0.8cm} \begin{cases}
\|\nabla^2  \partial_j\mathscr{B}  (\partial_i f)\|_{\bL^q(G)} \le 
c (\|\partial_i \nabla_\ast  \nabla  f\|_{\bL^q(G)} + \|\partial_j  \partial_n 
\partial_n f\|_{L^q(G)}
+ R^{-1}\|\nabla^2 f\|_{\bL^q(G)}) 
\\[0.3cm]
\forall\,  f\in W^{3,\, q}_0(G),
\end{cases}
\label{2.6-1}
\end{align}
$(i, j=1,\ldots,n)$ with a constant $c$, depending on $n, q$ and ${\rm ratio}(G)$ only.  Furthermore, 
from \eqref{2.3} we get 
\be\label{2.5-1}
\diver \mathscr{B} (f) (x)=f(x) -  \phi \Big(\frac {x} {R}\Big) R^{-n} \intl_{G } f(y) dy \quad \mbox{for a.\,e.} \,\,\, x\in G.
\ee

\section{Proof of Theorem\,1}
\label{sec:3}
\setcounter{secnum}{\value{section} \setcounter{equation}{0}
\renewcommand{\theequation}{\mbox{\arabic{secnum}.\arabic{equation}}}}

\vspace{0.5cm}
 {\bf Proof}  $1^\circ $  By decomposing the right-hand side into a solenoidal field, and 
a gradient field, we are able to reduce the problem to the case $\diver \bbf=0$. 
Let $\bE: \bW^{ k,\, q}(\Omega ) \rightarrow \bW^{ k,\, q}(\R^n)$   denote an extension operator such that  
\[
\|\bE\bv\|_{\bW^{k,\, q}(\R^n)}  \le
c \|\bv \|_{\bW^{k,\, q}(\Omega )}\quad \forall\,  \bv \in \bW^{ k,\, q}(\Omega )  .   
\]
Let  $\bP: \bW^{ k,\, q}(\R^n) \rightarrow \bW^{ k,\, q}_{0, \sigma }  (\R^n)$ denote the Helmholtz-Leray projection. Given $\bv\in \bW^{ k,\, q}(\Omega ) $
 we have 
\[
\bv = \bP\bE\bv + (I-\bP)\bE\bv \quad \mbox{a.\,e. in}\,\,\, \Omega. 
\]
In addition, there exists a constant $c>0$ depending only on $n, q,k$ and 
$\Omega $ such that 
\be\label{3.1}
\|\bP\bE\bv\|_{\bW^{ k,\, q}(\Omega )} + 
\|(I-\bP)\bE\bv\|_{\bW^{ k,\, q}(\Omega )}
\le  c \|\bv\|_{\bW^{ k,\, q}(\Omega )}\quad 
\forall\, \bv\in \bW^{ k,\, q}(\Omega ). 
\ee

\hspace*{0.5cm}
Now, for $ \bbf \in L^s(0,T; \bW^{k,\, r}(\Omega ))$ let  $(\bu, p)$ be a strong solution to \eqref{1.1}--\eqref{1.4}.  Observing $I-\bP = \nabla (\Delta ^{-1} \diver  )$ recalling the definition of $\bE$ we get 
\[
\bE\bbf = \bP\bE \bbf+ (I-\bP)\bE  \bbf=\bP\bE \bbf + \nabla (\Delta ^{-1} \diver \bE\bbf ) \quad \mbox{a.\,e. in}\,\,\, Q. 
\]
Since, $\nabla (\Delta ^{-1} \diver \bE\bbf ) =
(\Delta ^{-1} \nabla \diver \bE\bbf ) \in L^s(0,T;\bW^{ k,\, q}(\R^n))$  we see that 
$\bP\bE \bbf \in L^s(0,T;$ $\bW^{ k,\, q}(\R^n))$.  Thus, we can  replace $\bbf $ by the restriction of $\bP\bE \bbf $ 
on  $Q$, and $p$ by the restriction of 
$-\Delta ^{-1} \diver \bE\bbf +p$ on   $Q$.  
Hence, in what follows without loss of generality we may   
assume  that 
\be\label{3.2}
\diver \bbf =0,\quad\mbox{and}\quad\Delta p = 0 \quad \mbox{a.\,e. in }\,\,\, Q. 
\ee

$2^\circ$  Secondly, we recall a well-known result by Giga and Sohr \cite{GS}  which is 
the following

\vspace{0.5cm}
{\bf Lemma\,3.1}  {\it Let $\Omega = \R^n$, $\Omega = \R^n_+$, $\Omega $  bounded or $ \Omega $ an exterior domain with $\partial \Omega  \in C^2$. 
For every $\bg \in L^s(0,T; \bL^q_\sigma (\Omega ))$\, $(1<s,q, <+\infty)$   
there exists 
a unique solution $(\bv, \pi) \in L^s(0,T; \bW^{2,\, q}_{\rm loc} (\Omega ))\times 
L^s(0,T; W^{1,\, q}_{\rm loc} (\Omega ))$ to the Stokes problem 
\begin{align*}
\partial_t \bv - \Delta \bv &=- \nabla \pi  + \bg \quad  
\mbox{and}\quad \diver \bv =0 \quad \mbox{in}  \quad \Omega \times (0,T),
\\
\bv  &= 0 \quad\mbox{on}\quad \partial \Omega \times (0,T),
\\
\bv(0)  &=0 \quad \mbox{on}\quad \Omega \times \{0\},
\end{align*}
such that 
$\partial_t \bv, \partial_{i}\partial_j \bv, \nabla \pi \in L^s(0,T; \bL^q(\Omega )) $\, $(i,j=1,\ldots,n)$, and there holds,
\be\label{3.3}
\|\partial_t \bv \|_{L^s(0,T; \bL^q(\Omega ))}
+\|\nabla^2\bv\|_{L^s(0,T; \bL^q (\Omega ))} +
\|\nabla \pi  \|_{L^s(0,T; \bL^q(\Omega ))} \leq c\|\bg\|_
{L^s(0,T; \bL^q(\Omega ))},
\ee
where the constant  $c$ depends only on $n, s, q$ and $\Omega $. 
}

\vspace{0.2cm}
\hspace*{0.5cm} 
As a consequence of Lemma\,3.1 we get the existence of a unique solution 
$(\bu, p) \in  L^s(0,T; \bW^{2,\, q}_{\rm loc} (\Omega ))\times 
L^s(0,T; W^{1,\, q}_{\rm loc} (\Omega ))$  to the Stokes system \eqref{1.1}--\eqref{1.4}, 
such that 
  \be\label{3.3-1}
\|\partial_t \bu \|_{L^s(0,T; \bL^q(\Omega ))}
+\|\nabla^2\bu\|_{L^s(0,T; \bL^q (\Omega ))}+
\|\nabla p  \|_{L^s(0,T; \bL^q(\Omega ))} \le 
c\|\bbf\|_{L^s(0,T; \bL^q(\Omega ))}.
\ee

\vspace{0.5cm}
$3 ^\circ$ {\it Local estimates} 
We restrict ourself to case that $\Omega $ is an exterior domain. 
The opposite case can be treated in a similar way. 
Clearly, $G := \R^n\setminus \overline{\Omega} $ is a bounded domain. 
Let $G', G''$ are bounded open sets such that  $\overline{G} \subset G'$ and $\overline{G'}\subset G''$.  
Set $\Omega''= \R^3\setminus \overline{G''} $  and  
$\Omega' =  \R^3\setminus \overline{G'} $.  
Then, let  $\zeta \in C^\infty(\R^3)$   denote a cut-off function such that 
$\zeta \equiv 1$ on $\Omega''$, and $\zeta \equiv 0$ in $ G'$. 
In particular, $\supp(\nabla \zeta ) \subset G''\setminus G'$.  Observing 
$\diver (\bu(t) \zeta ) = \bu(t)\cdot \nabla \zeta$,   
it follows that $\supp (\bu(t)\cdot \nabla \zeta ) \subset\subset 
G''\setminus G'$ for a.\,e. $t\in (0,T)$. 

\hspace*{0.5cm}
Next, let $1<R<+\infty$ such that  $G'' \subset B_R$.  By 
$\mathscr{B}:  W^{k-1,\, q}_0 (B_R) \rightarrow 
\bW ^{ k,\, q}_0 (B_R)$  we denote the Bogowski\u{i} operator defined in Section\,2.  
We now define 
\[
\bz(t) =\mathscr{B} (\bu(t) \cdot \nabla \zeta),\quad t\in [0,T). 
\]
Let $t\in (0,t)$.  Since $\intl_{B_R} \bu(t) \cdot \nabla \zeta dx  =0$, in view of \eqref{2.5-1} we  have 
\[
\diver \bz(t) =\bu(t) \cdot \nabla \zeta \quad \mbox{a.\,e. in $B_R$}.
\]
Thanks to \eqref{2.6}, recalling that $ {\rm ratio}(B_R)=1$, there exists a constant $c>0$ depending only on $q$ and $n$ such that 
\[
\|\bz(t)\|_{\bW^{3,q}(B_R) } \le c 
\|\bu(t) \cdot \nabla \zeta \|_{W^{2,q}(B_R)} 
\quad \mbox{for a.\,e. $t\in (0,T)$}.
\]
Making use of the embedding $\bW^{3,q}_0(B_R) \hookrightarrow 
\bW^{3,q}(\R^n) $  the above inequality implies that 
$\bz \in L^s(0,T; \bW^{ 3,\, q}(\R^n))$.  Together with \eqref{3.3-1}, and  the 
Sobolev-Poincar\'{e} inequality we obtain  
\be\label{3.5}
\|\bz\|_{L^s(0,T; \bW^{3,q}(\R^n)) } \le c  
\|\bu \|_{L^s(0, T; \bW^{2,q}(\Omega \cap B_R))} \le 
c\|\bbf\|_{L^s(0,T; \bL^q(\Omega ))}.
\ee
By an analogous reasoning taking into account 
$\partial _t\bz = 
\mathscr{B} (\partial_t \bu \cdot \nabla \zeta) $  a.\,e. in $\R^n \times (0,T)$  we see that $\partial_t \bz\in L^s(0,T; \bW^{ 1,\, q}(\R^n))$. In addition, by virtue of  \eqref{3.3-1}  we obtain  
\be\label{3.6}
\|\partial_t \bz \|_{L^s(0,T; \bW^{1,q}(\R^3)) } \le c  
\|\partial_t \bu  \|_{L^s(0, T; \bL^{q}(\Omega ))} 
\le  c\|\bbf\|_{L^s(0,T; \bL^q(\Omega ))}.  
\ee

\hspace*{0.5cm}
Next, let  $k\in \{1,\ldots, n\}$ be fixed.  We define 
\[
\begin{cases}
\bv(x,t) = \partial_k (\bu(x,t) \zeta(x)  - \bz(x,t)), \quad & (x,t)\in 
(G''\setminus G')\times (0,T),  
\\[0.2cm]
\bv(x,t) =-\partial_k \bz(x,t),\quad   & (x,t)\in 
(\R^n \setminus (G''\setminus G'))\times (0,T),
\end{cases}
\]
and 
\[
\begin{cases}
\pi (x,t) =\partial_k (p(x,t) \zeta(x)), \quad & (x,t)\in 
(G''\setminus G')\times (0,T),  
\\[0.2cm]
\pi (x,t) =0,\quad   & (x,t)\in (\R^n \setminus (G''\setminus G'))\times (0,T).
\end{cases}
\]
Then the pair $(\bv, \pi ) $ solves the Stokes system 
\begin{align*}
\hspace*{3cm} \diver \bv &= 0  &&\mbox{in}\quad \R^n\times (0,T),
 \hspace*{4cm}
\\
\partial_t \bv - \Delta \bv &=- \nabla \pi  + \bg \quad  &&\mbox{in}\quad \R^n\times (0,T),
\\
\bv  &= 0 &&\mbox{on}\quad \R^n\times \{0\},
\end{align*}
where
\begin{align*}
\bg &= (p-p_{B_R} )\nabla \zeta - 2 \partial_k(\nabla \bu \cdot \nabla \zeta )- 
\partial_k(\bu \Delta \zeta)
\\
&\qquad \qquad- 
\partial_k\partial _t\bz + \partial_k\Delta \bz + \partial_k(\bbf \zeta )
\qquad \mbox{a.\,e. in $\R^n\times (0,T)$}.
\end{align*}
In view of \eqref{3.3}, \eqref{3.5}, and \eqref{3.6}  we 
see that $\bg \in L^s(0,T; \bL^q(\R^n))$. In addition, there holds
\[
\|\bg\|_{ L^s(0,T; \bL^q(\R^n))} \le  c 
\|\bbf\|_{ L^s(0,T; \bW^{ 1,\, q}(\Omega ))}.   
\]
Thus, applying Lemma\,3.1 with $\Omega = \R^n$, and using the last inequality we see that
\begin{align*}
&\|\partial_t \bv \|_{L^s(0,T; \bL^q(\R^n))}
+\|\nabla^2\bv\|_{L^s(0,T; \bL^{q}(\R^n))} +
\|\nabla \pi \|_{L^s(0,T; \bL^q(\R^n))} 
\\
&\qquad\le c\|\bg\|_{L^s(0,T; \bL^q(\R^n))}
\\
&\qquad \le c\|\bbf\|_{L^s(0,T; \bW^{ 1,\, q}(\Omega ))}. 
\end{align*}
Recalling the definition of $\bv$,  making use of \eqref{3.5}, \eqref{3.6}, and \eqref{3.3-1},  we infer from above 
\begin{align*}
&\Big\|\zeta \partial_t \partial_k\bu \|_{L^s(0,T; \bL^q(\Omega))}
+ \|\zeta  \nabla^2 \partial_k\bu\|_{L^s(0,T; \bL^{q}(\Omega ))} +
\|\zeta \nabla \partial_k p \|_{L^s(0,T; \bL^q(\Omega))} 
\\
&\qquad \le c\|\bbf\|_{L^s(0,T; \bW^{ 1,\, q}(\Omega ))}. 
\end{align*}

Iterating the above argument $k$ times, we get 
\begin{align}
& \|\partial_t \bu \|_{L^s(0,T; \bW^{k,\, q}(\Omega '))}
+ \| \bu\|_{L^s(0,T; \bW^{k+2,\, q}(\Omega'))} +
\|\nabla  p\|_{L^s(0,T; \bW^{k,\, q} (\Omega '))} 
\cr 
&\quad \le c\|\bbf\|_{L^s(0,T; \bW^{k,\, q}(\Omega ))}
\label{3.7}
\end{align}
$ (k\in \N)$, where $c=\const>0$, depending on $s, q, k$, and $\Omega $ only.

\vspace{0.3cm}
$4^\circ$ {\it Boundary regularity}  Let $x_0 \in \partial \Omega $. Up to translation and rotation we may assume that $x_0 =0$ and $\bn(0)= -\bfe_n$, where $\bn(0)$ denotes the outward unite normal on $\Omega $ at  $x_0$.  
According to our assumption 
on the boundary of $\Omega $  there exists $0<R< +\infty$,  and 
$h\in C^{2+k}(B_R') $ such that  
\begin{itemize}
\item[(i)] \quad $\partial \Omega \cap (B_R'\times (-R,R)) =
\{(y', h(y')); y'\in B_R'\};$

\item[(ii)] \quad $\{(y', y_n); y'\in B_R',  h(y') < y_n <h(y')+  R\} \subset \Omega; $

\item[(iii)] \quad $\{(y', y_n); y'\in B_R',   -R + h(y')< y_n < h(y')\} \subset \Omega^c $ \,\footnotemark\,.

\end{itemize}
\footnotetext{\,Here $y'=(y_1,\ldots, y_{n-1} )\in \R^{n-1} $, and $B_R'$  denotes the two dimensional ball $\{(y_1, \ldots, y_{n-1} ): y^2_1+\ldots + y_{n-1} ^2 < R^2\}$. }

Set $U_R = B_R'\times (-R,R), U_R^+ = B_R'\times (0,R)$,  and define 
$\bPhi: U_R \rightarrow \bPhi(U_R) $ by 
\[
\bPhi (y) = (y', h(y')+y_n)^{\top},\quad y\in U_R. 
\]
Elementary, 

\begin{align*}
D\bPhi(y)  &=\begin{pmatrix}
  1   & 0  &   \dots &  0 &  0 \\  
   0  &  1 &   \dots  & 0  &  0\\
  \vdots  & \vdots &   \ddots  & \vdots  &  \vdots \\
  & & & & \\
   0 & 0 &  \dots & 1& 0  \\
\partial_1 h(y)     &   \partial_2 h(y) &   \dots & \partial_{n-1}  h(y) &      1
\end{pmatrix}, 
\\[0.7cm]
(D\bPhi(y))^{-1}    &=  \begin{pmatrix}
  1   & 0  &   \dots &  0 &  0 \\  
   0  &  1 &   \dots  & 0  &  0\\
  \vdots  & \vdots &   \ddots  & \vdots  &  \vdots \\
  & & & & \\
   0 & 0 &  \dots & 1& 0  \\
-\partial_1 h(y)     &   -\partial_2 h(y) &   \dots & -\partial_{n-1}  h(y) &      1
\end{pmatrix}.
\end{align*}
For the outward unit  normal at  $x = \bPhi(y) $ we have 
\[
\bn(x)  = \bN(y) = 
\frac {(\partial_1h (y),\ldots,  \partial_{n-1} h (y), -1)} {\,\sqrt[]{1+ |\nabla h(y)|^2} }, \quad y \in  B_R'\times \{0\}.
\]

In addition, one calculates 
\be\label{3.8}
\partial_{x_i} \circ \bPhi =\partial_{y_i} - (\partial_{x_i} h)\partial_{y_n}  \quad \mbox{in}\quad U_R,\quad i=1,\ldots, n\,\footnotemark\,.
\ee
\footnotetext{\,Since $h$ is independent on $y_n$ there holds 
$\partial _{x_n} \circ \bPhi =\partial_{y_n}$.}

\hspace*{0.5cm}
We set $\bU =\bu \circ \bPhi, P = p  \circ \bPhi $ and 
$\bF = \bbf\circ \bPhi $   a.\,e.  in  
$U^+_R\times (0,T)$.  By the aid of \eqref{3.8}
we easily get  
\begin{align}
 ({\diver}_x \bu ) \circ \bPhi &= {\diver}_y \bU - 
\nabla h \cdot \partial_{y_n} \bU = 0,
\label{3.9}
\\
(\Delta_x \bu)\circ \bPhi &=\Delta_y \bU - 2\nabla h \cdot \nabla_y \partial_{y_n}\bU + |\nabla h|^2 \partial_{y_n} \partial_{y_n}\bU - 
(\Delta h) \partial_{y_n}\bU,
\label{3.10}
\\
(\nabla_x p) \circ \bPhi  &=
\nabla_y P - (\nabla h )\partial_{y_n} P,
\label{3.11}
\end{align}
a.\,e. in $U^+_R\times (0,T)$.    Firstly, owing to \eqref{3.9}  from the equation 
\eqref{1.1} we get 
\be\label{3.12}
{\diver}_y \bU =
\nabla h \cdot \partial_{y_n} \bU \quad \mbox{a.\,e. in}\quad U^+_R\times (0,T),
\ee
 and with help of \eqref{3.10} and 
\eqref{3.11}  the equation \eqref{1.2} turns into 
\begin{align}
\partial_t \bU - \Delta \bU 
%\cr 
%\\
%&\quad 
&=- \nabla P   +  (\partial_{y_n} P) 
 \nabla h- 2\nabla h \cdot \nabla \partial_{y_n}\bU + 
|\nabla h|^2 \partial_{y_n} \partial_{y_n}\bU 
\cr 
 &\qquad\qquad\qquad- (\Delta h)\partial_{y_n} \bU + \bF
\label{3.14}
\end{align}
a.\,e. in $U^+_R\times (0,T)$.  

\hspace*{0.5cm}
Note that the assumption $\bn(0)=-\bfe_n$   implies $\nabla h(0) =0$.  We now choose  
$0<\delta < +\infty$ sufficiently small, which will be specified later.  Since $\nabla h\in \bC^0(U_R)$, 
there exists $0< \rho  < \frac {R} {2}$ such that   
\be\label{3.14-1}
|\nabla h(y)| \le \delta \quad \forall\,  y\in U_{2\rho  }.    
\ee
Let $\zeta \in C^\infty_0(U_{2\rho } )$ denote a cut-off function such that  
$0\le \zeta \le 1$ in $U_{2\rho } $,  and $\zeta \equiv 1$ on $U_\rho $.  
We define $\widetilde{\bU}: \R^n_+ \times (0,T) \rightarrow \R^n$ by 
\[
\hspace*{-1cm} \widetilde{\bU}(y,t) = \zeta(y)\bU(y,t), \quad 
y\in U^+_{2\rho }\times (0,T), 
\quad \widetilde{\bU}(y,t) =0 \quad \mbox{if}\quad 
y\in \R^n_+\setminus U^+_{2\rho }\times (0,T).
\]
Let $\mathscr{B}: W^{ k-1,\, q}_0(U_{2\rho }^+ ) \rightarrow 
\bW^{ k,\, q}(\R^n_+)$  denote the Bogowski\u{i} operator defined in Section\,2. 
We set 
\begin{align*}
 \bz_1(y,t) \,&= \mathscr{B} (\zeta \nabla h\cdot \partial_{y_n} \bU)(y,t), 
\\
\bz_2(y,t) \,&=\mathscr{B} (\nabla \zeta\cdot \bU)(y,t),\quad 
(y, t) \in \R^n_+\times (0,T).
\end{align*}

Let $k\in \{1,\ldots, n-1\}$ be fixed.  We define 
\begin{align*}
\bV(y,t) &= \partial_{k}  (\widetilde{\bU}(y,t) - \bz_1(y,t) - \bz_2(y,t)),   
\\
\Pi(y,t) &=\partial_{k} (\zeta(y) P(y,t) ), 
\end{align*}
$(y, t) \in \R^n_+\times (0,T)$.  Observing that 
\[
\intl_{U^+_{2\rho } }   \zeta \nabla h\cdot 
\partial_{n} \bU(t) + \nabla \zeta  \cdot \bU(t) dy =\intl_{U^+_{2\rho } } {\diver}_{y}  \widetilde{\bU}(t) dy =0 
\quad \mbox{for a.\,e. $t\in (0,T)$},
\]
by the aid of  \eqref{2.5-1} we calculate 
\begin{align}
{\diver}_y \bV = \partial_{k}  \Big(\zeta \nabla h \cdot \partial_{n} \bU  + 
\nabla \zeta \cdot \bU - \zeta \nabla h\cdot \partial_{n} \bU - 
\nabla \zeta \cdot \bU \Big)   
%\cr 
%\,&
=0
\label{3.15}
\end{align}
a.\,e. in $\R^n_+\times (0,T)$.  In addition, taking into account \eqref{3.14}, we find 
\begin{align*}
\partial_t \bV - \Delta \bV   &= 
\partial_{k} \Big(\zeta \partial_t \bU - \zeta \Delta \bU - 2 \nabla \zeta \cdot \nabla \bU - (\Delta \zeta) \bU\Big)
\\
&\qquad - \partial_{k} (\partial_t \bz_1 - 
\Delta \bz_1\Big) - \partial_{k} 
(\partial_t \bz_2 - \Delta \bz_2)
\\
&=-\nabla  \Pi +\partial_{k} (\ (P-P_{U^+_{2\rho }}) \nabla \zeta )  -  
\partial_{k} \Big(2 \nabla \zeta \cdot \nabla \bU + (\Delta \zeta) \bU\Big)
\\
&\qquad -\partial_{k} (\partial_t \bz_1 - 
\Delta \bz_1) -\partial_{k} 
(\partial_t \bz_2 - \Delta \bz_2)
\\
&\qquad +
\partial_{k} \Big(\zeta  (\partial_{n} P) 
  \nabla h- 2\zeta \nabla h \cdot \nabla \partial_{n}\bU + 
\zeta |\nabla h|^2 \partial_{n} \partial_{n}\bU 
\\
&\qquad \qquad \qquad  - \zeta  (\Delta h)\partial_{n} \bU + \zeta \bF\Big).
\end{align*}
Thus, $(\bV, \Pi )$  solves the following Stokes system 
\begin{align*}
\hspace*{2cm} \diver \bV &= 0 \quad &&\mbox{in}\quad \R^n_+\times (0,T), \hspace*{4cm}
\\
\partial_t \bV - \Delta \bV &= 
-\nabla \Pi + \bG  &&\mbox{in}\quad \R^n_+\times (0,T),
\\
\bV \,&= 0 \,\, &&\mbox{on}\quad \partial \R^n_+\times (0,T),
\end{align*}
where $\bG = \bG_1 + \ldots +\bG_6$ with 
\begin{align*}
\bG_1 &= \partial_{y_k} ((P-P_{U_{2\rho }^+ })  \nabla \zeta),     
\\
\bG_2 &=-
\partial_{k} (2 \nabla \zeta \cdot \nabla \bU + (\Delta \zeta) \bU),
\\
\bG_3 &= -\partial_{k} (\partial_t \bz_1- \Delta \bz_1\Big),
\\
\bG_4 &=-\partial_{k} (\partial_t \bz_2 - \Delta \bz_2),
\\
\bG_5 &=\partial_{k} \Big(\zeta 
 (\partial_{n} P) 
  \nabla h- 2\zeta \nabla h \cdot \nabla \partial_{n}\bU + 
\zeta |\nabla h|^2 \partial_{n} \partial_{n}\bU \Big),
\\
\bG_6 &=  \partial_{k} (- \zeta  (\Delta h)\partial_{n} \bU 
+ \zeta \bF).
\end{align*}

\hspace*{0.5cm} 
In what follows we shall establish some important  estimates of $\bz_1$ and $\bz_2$, where we will make essential use 
of the properties of $\mathscr{B}$ (cf. Section\,2).   Starting with $\bz_1$, we write $\bz_1 = \bz_{1,1}+ \bz_{1,2}$, where 
\[
\bz_{1,1} =\mathscr{B} (\partial_{n} (\zeta \nabla h\cdot \bU)),\quad
\bz_{1,2} =-\mathscr{B} ((\partial_{n} \zeta) \nabla h\cdot \bU)).
\]
Let $t\in (0,T)$ be fixed.  Using \eqref{2.5}, \eqref{2.6} with $j=k, i=n$  and $f=\zeta \nabla h\cdot \bU $, and observing 
 $\partial_t  \mathscr{B} = 
\mathscr{B} \partial _t$,  we see that  
\begin{align*}
 \|\partial_t \partial_k\bz_1 (t)\|_{\bL^q(\R^n_+)} 
\,&\leq\,   \|\partial_t \partial_k\bz_{1,1} (t) \|_{\bL^q(\R^n_+)}
+\|\partial_t \partial_k\bz_{1,2} (t) \|_{\bL^q(\R^n_+)}
\\
& \le 
c \|\partial_t\partial _k (\zeta \nabla h\cdot \bU) (t) \|_{L^q(\R^n_+)} 
+
c\rho ^{-1} \|\partial_t \bU (t)\|_{\bL^q(\R^n_+)}
\\
& \le  c \delta \|\partial_t  
\partial_k \widetilde{\bU}  (t)\|_{L^q(\R^n_+)} +
c(\|h\|_{C^2} + \rho ^{-1})  \|\partial_t \bU (t) \|_{\bL^q(U_{2\rho }^+ )}.
\end{align*}
Taking the above inequality to   the $s$-th power,  and  integrating the resulting equation  in time 
over $(0,T)$,  we get 
\begin{align}
& \|\partial_t \partial_k\bz_1 \|_{L^s(0,T; \bL^q(\R^n_+))}
\cr 
&\quad \le c\delta \|\partial_t \partial_k 
\widetilde{\bU} \|_{L^s(0,T; \bL^q(\R^n_+))} 
+ c(\|h\|_{C^2} + \rho ^{-1})  \|\partial_t \bU \|_{L^s(0,T; \bL^q(U_R^+))}.
\label{3.16}
\end{align}
On the other hand,  using \eqref{2.5}, \eqref{2.6-1} with $j=k, i=n$,  and 
$f=\zeta \nabla h\cdot \bU(t) $,  we see that 
\begin{align*}
\hspace*{-1cm} \|\nabla^2 \partial_k\bz_1(t)\|_{\bL^q(\R^n_+)}   \,&\leq\, 
c \|\partial_n \nabla_\ast \nabla (\zeta  \nabla h \cdot \bU)(t) \|_{L^q(\R^n_+)} 
+ c \|\partial_n \partial_n \partial_k (\zeta  \nabla h \cdot \bU)(t)  \|_{L^q(\R^n_+)}
\\
&\quad + c \rho ^{-1} \|\nabla^2(\zeta  \nabla h \cdot \bU)(t) \|_{L^q(\R^n_+)}  
+ c\| \nabla^2 ((\partial_n\zeta) \nabla h \cdot \bU)(t)\|_{L^q(\R^n_+)}.
\end{align*}
By means of product rule and Poincar\'{e}'s inequality we find 
\begin{align*}
\|\nabla^2 \partial_k\bz_1(t)\|_{\bL^q(\R^n_+)} &\leq
c \delta \|\nabla^2 \nabla_\ast \widetilde{\bU}(t) \|_{L^q(\R^n_+)}
+ c(\|h\|_{C^3} + \rho ^{-1}) \|\nabla^2\bU(t)\|_{\bL^q(U_R^+)}.
\end{align*}
We now take the above inequality to the $s$-th power, integrating the result in time  over $(0,T)$, we obtain  
\begin{align}
\|\nabla^2 \partial_k\bz_1\|_{L^s(0,T; \bL^q(\R^n_+))} \,&\leq    
c \delta \|\nabla^2 \nabla_\ast \widetilde{\bU} \|_{ L^s(0,T; L^q(\R^n_+))}
\cr 
&\qquad + c(\|h\|_{C^3} + \rho ^{-1}) \|\nabla^2\bU\|_{ L^s(0,T; \bL^q(U_R^+))}.
\label{3.17}
\end{align}

\hspace*{0.5cm} 
By an analogous reasoning, making use of \eqref{2.5}, and Poincare's inequality, we infer 
\begin{align}
&\|\partial_t \bz_2  \|_{L^s(0,T; \bL^q(\R^n_+))} 
+\|\nabla^2 \bz_2\|_{L^s(0,T; \bL^q(\R^n_+))}   
\cr 
&\qquad \qquad  \leq  c \rho^{-1}  
\Big(\|\partial_t \bU \|_{L^s(0,T; \bL^q(U_R^+))} 
+\|\nabla^2 \bU\|_{L^s(0,T; \bL^q(U_R^+))} \Big).
\label{3.18}
\end{align}

\hspace*{0.5cm} 
We are now in a position to estimate $\bG_1, \ldots, \bG_6$.  First by virtue of Poincar\'{e}'s inequality we easily estimate  
\[
\|\bG_1\|_{L^s(0,T; \bL^q(\R^n_+))} \le c \rho ^{-1} 
\|\nabla P\|_{L^s(0,T; \bL^q(U_R^+))}.   
\]
Analogously, 
\[
\|\bG_2\|_{L^s(0,T; \bL^q(\R^n_+))} \le  c \rho ^{-1} 
\|\nabla^2\bU\|_{L^s(0,T; \bL^q(U_R^+))}.   
\]

\hspace{0.5cm}
Next,  with the help of \eqref{3.16}, \eqref{3.17}, and \eqref{3.18}  we see that 
\begin{align*}
&\|\bG_{3}  \|_{L^s(0,T; \bL^q(\R^n_+))} +
\|\bG_{4}  \|_{L^s(0,T; \bL^q(\R^n_+))}
\\  
&\quad \le c \delta  \Big(\|\partial_t \partial_k\widetilde{\bU} \|_
{L^s(0,T; \bL^q(\R^n_+))} + \|\nabla^2 \nabla_\ast \widetilde{\bU}\|_
{L^s(0,T; \bL^q(\R^n_+))}\Big)
\\
&\qquad + c (\|h\|_{C^3}+ \rho ^{-1}  ) \Big(\|\partial_t \partial_k\widetilde{\bU} \|_{L^s(0,T; \bL^q(U_R^+))} + 
\|\nabla^2 \nabla_\ast \widetilde{\bU}\|_{L^s(0,T; \bL^q(U_R^+))}\Big).
\end{align*}

\hspace{0.5cm}
Then applying the product rule, and using Poincar\'{e}'s inequality, we get 
\begin{align*}
\|\bG_5 \|_{L^s(0,T;\bL^q(\R^n_+))} &\le 
 c \delta \Big(\|\nabla \Pi \|_{L^s(0,T;\bL^q(\R^n_+))}+
 \|\nabla_\ast \nabla^2 \widetilde{\bU}\| _{L^s(0,T;\bL^q(\R^n_+))} \Big)
\\
&\quad + c(\|h\|_{C^2} + \rho ^{-1} ) \Big(\|\nabla P\|_{L^s(0,T;\bL^q(U_R^+))}  
+ \|\nabla^2\bU\|_{L^s(0,T;\bL^q(U_R^+))}\Big).
\end{align*}  

\hspace*{0.5cm} 
Finally,  we estimate 
\begin{align*}
\|\bG_6 \|_{L^s(0,T;\bL^q(\R^n_+))} &\leq 
c(\|h\|_{C^3} + \rho ^{-1} ) \Big(\|\nabla P\|_{L^s(0,T;\bL^q(U_R^+))}  
+\|\nabla^2\bU\|_{L^s(0,T;\bL^q(U_R^+))}\Big).  
\\
&\qquad +c \rho ^{-1} \|\bF\|_{L^s(0,T; \bL ^q(U_R^+))} +c\|\partial_k\bF\|_{L^s(0,T; \bL ^q(U_R^+))}.
\end{align*}

Appealing to  Lemma\,3.1 (cf. \cite{GS}) for the case $\Omega = \R^n_+$  using the above estimates for  
$\bG_1, \ldots, \bG _6$, we obtain
\begin{align*}
 & \|\partial_t \bV \|_{L^s(0,T; \bL^q(\R^3_+))} 
+\|\nabla^2 \bV\|_{L^s(0,T; \bL^q(\R^3_+))}   +  
\| \nabla \Pi \|_{L^s(0,T; \bL^q(\R^3_+))} 
\\
&\quad \le c \|\bG_1 + \ldots +\bG_6\|_{L^s(0,T;\bL^q(\R^n_+))}
\\
&\quad \le c \delta 
\Big(\|\partial_t \nabla_\ast \widetilde{\bU} \|_{L^s(0,T; \bL^q(\R^3_+))} 
+\|\nabla^2 \nabla_\ast \widetilde{\bU}\|_{L^s(0,T; \bL^q(\R^3_+))}   +  
\| \nabla \Pi \|_{L^s(0,T; \bL^q(\R^3_+))}\Big)
\\
&\qquad +  c(\|h\|_{C^3} + \rho ^{-1} ) 
\Big(\|\partial_t \bU \|_{L^s(0,T;\bL^q(U^+_R))} 
+\|\nabla^2 \bU\|_{L^s(0,T;\bL^q(U^+_R))}
\\
&\qquad \qquad\qquad\qquad\qquad +   \|\nabla P\|_{L^s(0,T;\bL^q(U^+_R))} +\|\bF\|_{L^s(0,T; \bW^{ 1,\, q}(U_R^+))}\Big).
\end{align*}
Recalling $\bV = \partial_k(\widetilde{\bU}- \bz_1-\bz_2)$,   making use of \eqref{3.16}, 
\eqref{3.17} and \eqref{3.18},  from the last inequality we infer 
\begin{align}
 & \|\partial_t \nabla_\ast \widetilde{\bU} \|_{L^s(0,T; \bL^q(\R^3_+))} 
+\|\nabla^2 \nabla_\ast\widetilde{\bU}\|_{L^s(0,T; \bL^q(\R^3_+))}   +  
\| \nabla \Pi \|_{L^s(0,T; \bL^q(\R^3_+))} 
\cr 
&\quad \le c_0 \delta 
\Big(\|\partial_t \nabla_\ast \widetilde{\bU} \|_{L^s(0,T; \bL^q(\R^3_+))} 
+\|\nabla^2 \nabla_\ast \widetilde{\bU}\|_{L^s(0,T; \bL^q(\R^3_+))}  
+\| \nabla \Pi \|_{L^s(0,T; \bL^q(\R^3_+))}\Big)
\cr 
&\qquad + c_1 
\Big(\|\partial_t \bU \|_{L^s(0,T;\bL^q(U^+_R))} 
+\|\nabla^2 \bU\|_{L^s(0,T;\bL^q(U^+_R))}
\cr 
&\qquad \qquad\qquad\qquad\qquad +  \|\nabla P\|_{L^s(0,T;\bL^q(U^+_R))} +\|\bF\|_{L^s(0,T; \bW^{ 1,\, q}(U_R^+))}\Big),
\label{3.19}
\end{align}
where $c_0=c_0(n,q, s)$ and $c_1=c_1(n,q,s, \|h\|_{C^3}, \rho )$.  
On the other hand, recalling the definition of $\bU, P$, and $\bF$, with the help of   
\eqref{3.10}, \eqref{3.11},  and \eqref{3.7} we find 
\begin{align}
& \|\partial_t \bU \|_{L^s(0,T;\bL^q(U^+_R))} 
+\|\nabla^2 \bU\|_{L^s(0,T;\bL^q(U^+_R))}  
\cr 
&\qquad \qquad +  c \|\nabla P\|_{L^s(0,T;\bL^q(U^+_R))} +
\|\bF\|_{L^s(0,T; \bW^{ 1,\, q}(U_R^+))}  \leq  c 
\|\bbf\|_{L^s(0,T; \bW^{ 1,\, q}(\Omega ))} 
\label{3.20}
\end{align}
with a constant $c$ depending on $n, q, s$ and $h$.  
Now, in \eqref{3.19} we take   $\delta = \frac {1} {2c_0}$ and estimate  the right-hand side of \eqref{3.19} by the aid of   \eqref{3.20}. 
This leads to   
\begin{align*}
 & \|\partial_t \nabla_\ast \widetilde{\bU} \|_{L^s(0,T; \bL^q(\R^3_+))} 
+\|\nabla^2 \nabla_\ast\widetilde{\bU}\|_{L^s(0,T; \bL^q(\R^3_+))}   +  
\| \nabla \Pi \|_{L^s(0,T; \bL^q(\R^3_+))} 
\\
&\qquad \le  c_2\| \bbf\|_{L^s(0,T;\bW^{ 1,\, q}(\Omega ))},
\end{align*}
where $c_2=c_2(n,q,s, \|h\|_{C^3}, \rho)$. 

\hspace*{0.5cm} 
By a standard iteration argument we obtain 
\begin{align}
 & \|\partial_t \nabla_\ast^k \bU \|_{L^s(0,T; \bL^q(U_\rho ^+))} 
+\|\nabla^2 \nabla_\ast^k\bU\|_{L^s(0,T; \bL^q(U_\rho ^+))}  
+\| \nabla \nabla_\ast^k P\|_{L^s(0,T; \bL^q(U_\rho ^+))} 
\cr 
&\qquad \le c \| \bbf\|_{L^s(0,T;\bW^{ 1,\, q}(\Omega ))},
\label{3.21}
\end{align}
where $c=\const$ depending only on $n, q, s, k, \|h\|_{C^{k+2} }$ and $\rho$.

\vspace{0.2cm}
$5^\circ$  {\it Estimation of the full pressure gradient}  Recalling that $\Delta_x p = 0 $,  
with the help of \eqref{3.10} we calculate 
\begin{align*}
0 =\Delta_x p \circ \bPhi &=\Delta_y P- 2\nabla h \cdot 
\nabla \partial_{y_n}P  + 
|\nabla h|^2 \partial_{y_n} \partial_{y_n}P- (\Delta h )\partial_{y_n}P  
\\
&=(1+ |\nabla h|^2) \partial_{n} \partial_{y_n}P  + \Delta'_y P- 2\nabla h \cdot 
\nabla_\ast \partial_{n}P  - (\Delta h )\partial_{y_n}P\,\footnotemark\,
\end{align*}
a.\,e. in  $ U^+_R $. Thus, 
\[
  (1+ |\nabla h|^2) \partial_{y_n} \partial_{y_n}P
=- \Delta '_y P  + 2\nabla h \cdot \nabla_\ast \partial_{y_n}P 
+ (\Delta h) \partial_{y_n}P 
\]
\footnotetext{\,Here $\Delta '_y$ stands for the differential operator $\partial_{y_1} \partial_{y_1}  + \ldots + 
\partial_{y_{ n-1} }\partial_{y_{ n-1}} $. }
a.\,e. in  $ U^+_R $. From this identity along with \eqref{3.21} with $k=1$  
it follows that 
\begin{align*}
\| \nabla^2_y P \|_{L^s(0,T; \bL^q (U_{\rho }^+ ))} &\leq 
c\Big(\|\nabla_{\ast}  \partial_{y_n}  P \|_{L^s(0,T; \bL^q (U_{\rho }^+  ))} 
+\|\nabla_y  P \|_{L^s(0,T; \bL^q (U_{\rho }^+))}\Big) 
\\
&\leq
c \|\bbf\|_{L^s(0,T; \bW^{1, q} (\Omega ))}.  
\end{align*}

Choosing $ \rho \in  \Big(0, \frac{R}{2}\Big)$ sufficiently small, and applying the above argument $k$-times, we get
\be\label{3.22}
\| \nabla^{k+1}_y  P\|_{L^s(0,T; \bL^q (U_{\rho }^+))} 
\le c \|\bbf\|_{L^s(0,T; \bW^{k, q} (\Omega ))}
\ee
with a constant $c$ depending on $n, q, s, k, \|h\|_{C^{k+2} } $, and $\rho$.  

\hspace*{0.5cm} 
Finally a standard covering argument, together with  \eqref{3.22}, and \eqref{3.7}  gives the estimate \eqref{1.5}, which completes 
the  proof of the Theorem\,1.\hfill \Beweisende

\vspace{0.5cm}
 \vspace{0.5cm}
{\bf Acknowledgements}
The present research has been supported by the German Research Foundation (DFG) through the project WO1988/1-1; 612414.

\end{document}